\documentclass[11pt,twoside]{amsart}

\usepackage{amsfonts}

\title[derived equivalences and Reynolds ideals]{Invariance
of generalised Reynolds ideals\\ under derived equivalences}

\author[Alexander Zimmermann]{Alexander Zimmermann}
\address{Universit\'e de Picardie,\newline  Facult\'e de Math\'ematiques et
LAMFA (UMR 6140 du CNRS),\newline    33 rue St Leu, \newline  F-80039
Amiens Cedex 1, \newline  France}
\email{alexander.zimmermann@u-picardie.fr}
\date{February 17, 2005; revised: August 25, 2005}

\newtheorem{Theo1}{{Theorem}}
\newtheorem*{Theo2}{{Theorem}}
\newtheorem{Lemma1}{{Lemma}}
\newtheorem{Def1}{{Definition}}
\newtheorem{Prop1}[Lemma1]{{Proposition}}
\newtheorem{Claim1}[Lemma1]{{Claim}}
\newtheorem{Rem1}{{Remark}}[section]
\newtheorem{Cor1}[Lemma1]{{Corollary}}
\newtheorem{Ex1}[Rem1]{{Example}}

\newenvironment{Rem}{\begin{Rem1}\rm}{\end{Rem1}}
\newenvironment{Theorem}{\begin{Theo1}}{\end{Theo1}}

\newenvironment{Claim}{\begin{Claim1}}{\end{Claim1}}

\newcommand{\uar}{\uparrow}
\newcommand{\dar}{\downarrow}
\newcommand{\lra}{\longrightarrow}

\newcommand{\ra}{\rightarrow}
\newcommand{\sdp}{\times\kern-.2em\vrule height1.1ex depth-.05ex}
\newcommand{\epi}{\lra \kern-.8em\ra}

\newcommand{\Z}{{\mathbb Z}}

\newcommand{\dickebox}{{\vrule height5pt width5pt depth0pt}}

 \normalbaselines
\begin{document}

\begin{abstract}
For any algebraically closed field $k$ of
positive characteristic $p$ and any non negative integer
$n$ K\"ulshammer defined ideals $T_nA^\perp$ of the centre of a
symmetric $k$-algebra $A$. We show that for derived equivalent
algebras $A$ and $B$
there is an isomorphism of the centres of $A$
and $B$ mapping $T_nA^\perp$ to $T_nB^\perp$ for all $n$.
Recently H\'ethelyi, Horv\'ath, K\"ulshammer and
Murray showed that this holds for Morita equivalent algebras.
\end{abstract}

\maketitle

\section*{Introduction}

Let $k$ be an algebraically closed field of characteristic $p>0$
and let $A$ be a finite dimensional symmetric $k$-algebra with non
degenerate symmetrising bilinear form $(\;,\;)$ on $A$.
K\"ulshammer defined in \cite{Kueldeutsch} ideals $T_nA^{\perp}$
of the centre of $A$ by the following construction. Let $KA$ be
the $k$-subspace of $A$ generated by $ab-ba$ for all $a,b\in A$
and set $T_nA:=\{x\in A\;|\; x^{p^n}\in KA\}$. Let $T_nA^\perp$ is
the subspace orthogonal to $T_nA$ with respect to the form
$(\;,\;)$ on $A$. Note that $T_nA^\perp$ is then an ideal of $ZA$
as $ZA=KA^\perp$ and that $T_nA$ is a $ZA$ submodule of $A$

In \cite{Kuelprog}  K\"ulshammer  shows that the equation
$(\zeta_n(z),x)^{p^n}=(z,x^{p^n})$ for any $x,z$ in the centre of
$A$ defines a mapping $\zeta_n$ from the centre of $A$ to the
centre of $A$. Moreover, $\zeta_n(A)=T_nA^{\perp}$. Many properties of
group algebras can be shown using the ideals $T_nA^\perp$.
Concerning the ideals
$T_nA^{\perp}$, H\'ethelyi et al. show in \cite{Kuel} that
$Z_0A\subseteq (T_1A^{\perp})^2\subseteq HA$, where $HA$ is the
Higman ideal of $A$, and where $Z_0A$ is the sum of the centres of
those blocks of $A$ which are simple algebras. They show that for
odd $p$ the left inclusion is an equality, whereas for $p=2$ one
gets $Z_0A=(T_1A^{\perp})^3=(T_1A^{\perp})\cdot (T_2A^{\perp})$.
Finally, the authors show that $e\cdot (T_nA^\perp)\cdot
e=T_n(eAe)^\perp$ for any idempotent $e$ of $A$. Now, the authors
use the fact that within all algebras Morita equivalent to $A$
there is an, up to isomorphism unique, smallest algebra $B=eAe$
Morita equivalent to $A$, the basic algebra. If $A$ is symmetric,
$B$ is symmetric as well, as follows in a more general context by
\cite{rogquest}. Multiplication by this idempotent induces an
isomorphism between the centers of an algebra $A$ and its basic
algebra $B$. Hence, the corresponding ideals $T_nA^{\perp}$ and
$T_nB^{\perp}$ are sent to each other by this isomorphism.
Composing two of them gives a corresponding statement for Morita
equivalent algebras.

In \cite[question 5.4]{Kuel} H\'ethelyi et al. ask whether for two
symmetric algebras $A$ and $B$, the condition that the derived
categories of $A$ and of $B$ are equivalent imply the existence of
an isomorphism $\varphi$ of their centres so that $\varphi$
induces an isomorphism between the ideals $T_nA^{\perp}$ and
$T_nB^{\perp}$ for all $n\in{\mathbb N}$. The main objective of
this paper is to give a positive answer to this question.

This way we provide new invariants for an equivalences between the
triangulated categories $D^b(A)$ and $D^b(B)$ for algebras $A$ and
$B$. A number of invariants are known. Suppose $D^b(A)\simeq
D^b(B)$ as triangulated categories, then we get an isomorphism of
the Hochschild homology $HH_*(A)\simeq HH_*(B)$ and the
Hochschild cohomology $HH^*(A)\simeq HH^*(B)$ (cf Rickard \cite{Ri3}),
the cyclic homology $HC_*(A)\simeq HC_*(B)$, the
cyclic cohomology $HC^*(A)\simeq HC^*(B)$ of the algebra $A$
(Keller~\cite{kellercyclic}), or by a result of Thomason and
Trobaugh the $K$-theory $K_*(A)\simeq K_*(B)$. Some of them are
quite useful in specializing the degree. So is $HH^0(A)\simeq
Z(A)$ the centre of the algebra $A$, or $rank_\Z(K_0(A))$ equals
the number of isomorphism classes of simple $A$-modules.
Nevertheless, if these few computable invariants coincide, it is
in general very difficult to decide whether two algebras have
equivalent derived categories or not. So, invariants which are
more easy to determine in examples will be very welcome. Our
result provides some of them.

The main result Theorem~\ref{main} will be proven in
Section~\ref{derivedinvariance}. Since there is no analogue of a
basic algebra for derived equivalences, we need to proceed
differently from H\'ethelyi's et al.'s proof for Morita
equivalence. Section~\ref{sec0} recalls some of the relevant
notation and results from homological algebra, for the convenience
of the reader. We use the characterisation \cite[(46)]{Kuelprog};
or \cite[Lemma 2.1]{Kuel} of $T_nA^\perp$ as the image of the
mapping $\zeta_n^A$ and define in Section \ref{sec1} the mapping
$\zeta_n^A$ in a functorial manner by means of a composition of
mappings between $A\otimes_kA^{op}$-modules. We apply the derived
equivalence to each of the factors and using results in
\cite{rogquest}, and some delicate commutativity considerations we
are able to show that the mapping induced by a standard derived
equivalence on the morphism sets are indeed as asked. For
notations concerning derived categories and equivalences we follow
\cite{derbuch}. Other references covering the needed background
are for example Gelfand-Manin~\cite{GelfandManin} or
Weibel~\cite{Weibel}. The notation may differ slightly there.

\medskip

\noindent
{\bf Acknowledgement:}
I want to thank Burkhard K\"ulshammer for having had a look into a first
draft of the manuscript, for pointing out a mistake there and for most
helpful suggestions how to fix this gap. I want to thank
John Murray for explaining to me during the AMS/DMV-meeting in May 2005
some of the historical background
around the origins of the ideals $T_nA^\perp$.

\section{A crash course on the relevant homological algebra}
\label{sec0}\label{somenotation}
For the reader's convenience and to fix notation
we shall recall some basic facts in
homological algebra as it is needed in the sequel. Basic source is
the book \cite{derbuch}, and for some more general aspects
Gelfand-Manin~\cite{GelfandManin}, Weibel~\cite{Weibel},
or Rickard~\cite{Ri3} as well as \cite{rogquest}.

For a commutative noetherian ring $k$ and a finitely generated
$k$-algebra $A$ we denote the category of finitely generated left
$A$-modules by $A-mod$ and the category of all $A$-modules
by $A-Mod$. Let $K(A-mod)$ be the category of complexes
in $A-mod$ modulo homotopy.
Recall that the derived category $D^b(A)$
of bounded complexes of finitely
generated $A$-modules is formed by bounded complexes in $K(A-mod)$
and formally inverting morphisms which induce isomorphisms on homology.
Recall furthermore that $A-mod$ is a full subcategory of $D^b(A)$
by mapping a module $M$ to a complex with homogeneous components $0$
in all degrees except in degree $0$ where the homogeneous component
is $M$  (cf e.g. Gelfand-Manin \cite[III \S 5 Proposition 2]{GelfandManin}).
Hence, any two objects $M$ and $N$ of $A$-mod may be considered as
object in $D^b(A)$, and then $Hom_{D^b(A)}(M,N)=Hom_A(M,N)$.
This fact will be used at various places.

Recall that in case $X$ is a complex in $D^b(A)$ whose homogeneous
components are all projective, then for any complex $Y$ in $D^b(A^{op})$
one has $Y\otimes_AX=Y\otimes_A^{\mathbb L}X$. In case $A$ and $B$ are
two algebras over a field $k$, then for any $X$ in
$D^b(A\otimes_kB^{op})$ there is an $\tilde X$ in $D^b(A\otimes_kB^{op})$
so that $X\simeq \tilde X$ and so that all homogeneous components
of $\tilde X$ are projective as $A$-modules and as $B^{op}$-modules
(cf \cite[Lemma 6.3.12]{derbuch}).

Let $B$ be a $k$-algebra which is projective as $k$-module.
By a result due to Keller (cf e.g. \cite{bernhard}) or
\cite[Chapter 8]{derbuch}) $D^b(A)$ is
equivalent to $D^b(B)$ as triangulated categories if and only if there is
a complex $X$ in $D^b(B\otimes_kA)$ so that
$X\otimes_A^{\mathbb L}-:D^b(A)\lra D^b(B)$ is an equivalence.
Such equivalences are called standard and $X$ is called
(two-sided) tilting complex.
If $B$ is symmetric, then $A$ is symmetric as well (cf \cite{rogquest})
and then the inverse equivalence to
$X\otimes_A^{\mathbb L}-$ is given by $Hom_k(X,k)\otimes_B^{\mathbb L}-$.
Moreover, Rickard has shown that this
$X\otimes_A^{\mathbb L}-\otimes_A^{\mathbb L}Hom_k(X,k)$
actually defines an equivalence
$D^b(A\otimes_kA^{op})\lra D^b(B\otimes_kB^{op})$ where
the $(A\otimes_kA^{op})$-module $A$ is mapped to $B$ (cf Rickard \cite{Ri3};
or \cite[Proposition 6.2.6]{derbuch}).
Hence $X$ induces an isomorphism
\begin{eqnarray*}
Z(A)=End_{A\otimes_kA^{op}}(A)&=End_{D^b(A\otimes_kA^{op})}(A)&
\simeq\\& End_{D^b(B\otimes_kB^{op})}(B)=&End_{B\otimes_kB^{op}}(B)=Z(B).
\end{eqnarray*}
This isomorphism is explicitly exhibited in
\cite[Proposition 6.2.6]{derbuch}.
In \cite{rogquest} it is shown that under the equivalence
induced by tensoring with $X$
the $(A\otimes_kA^{op})$-module $Hom_k(A,k)$ is mapped to $Hom_k(B,k)$.

We finish with some notation.
Let $\mathcal C$ be a category
and let  $X$, $Y$ and $Z$ be any three objects in $\mathcal C$.
We denote for any morphism $\varphi\in Hom_{\mathcal C}(X,Y)$
the induced mapping
$Hom_{\mathcal C}(Z,\varphi):Hom_{\mathcal C}(Z,X)\lra
Hom_{\mathcal C}(Z,Y)$ which is defined by
$\left(Hom_{\mathcal C}(Z,\varphi)\right)(\psi):=\varphi\circ\psi$
for any $\psi\in Hom_{\mathcal C}(Z,X)$. If $Z$ is clear from the context,
we write $Hom_{\mathcal C}(Z,\varphi)=:\varphi_*$ for short.

\section{Interpreting  $\zeta$}
\label{sec1}

Recall from Section~\ref{sec0} that
$Hom_{D^b(A\otimes_kA^{op})}(A,A)=
End_{A\otimes_kA^{op}}(A)\simeq Z(A)\;.$

Furthermore, by the adjointness formulas (cf e.g. Mac Lane
\cite[VI (8.7)]{Maclane}), we get
\begin{eqnarray*}
Hom_k(A\otimes_{A\otimes_kA^{op}}A,k)
&\simeq&Hom_{A\otimes_kA^{op}}(A,Hom_k(A,k))\\
f&\mapsto&\left(a\mapsto \left(b\mapsto f(a\otimes b) \right)\right)
\end{eqnarray*}
and since canonically by the very definition of a tensor product
$A\otimes_{A\otimes_kA^{op}}A\simeq A/KA$
where $KA=\sum_{a,b\in A}k\cdot (ab-ba)$ is the $k$-vector space
generated by commutators, we have a functorial isomorphism
\begin{eqnarray*}
Hom_k(A/KA,k) &\simeq&Hom_{A\otimes_kA^{op}}(A,Hom_k(A,k))\\
f&\mapsto&(a\mapsto(b\mapsto f(a b) ) )
\end{eqnarray*}

The mapping $A/KA\ni a\mapsto a^p\in A/KA$ was first defined by
Richard Brauer who called it the Frobenius mapping and proved that it
is well defined (cf K\"uls\-hammer \cite[II]{Kueldeutsch}) and semilinear.
Denote by $k^{(n)}$ the $n$ times Frobenius twisted copy of $k$.

The Frobenius mapping induces a well defined mapping
\begin{eqnarray*}
Hom_k(A/KA,k)&\lra&Hom_k(A/KA,k^{(1)})\\
f&\mapsto&\left(a\mapsto f(a^p)\right)
\end{eqnarray*}
The mapping
\begin{eqnarray*}
{Fr^k}_*:\;Hom_k(A/KA,k)&\lra&Hom_k(A/KA,k^{(1)})\\
f&\mapsto&\left(a\mapsto f(a)^p\right)
\end{eqnarray*}
induces a mapping
$$Hom_{ A\otimes_kA^{op}}(A,Hom_k(A,k))\lra
Hom_{ A\otimes_kA^{op}}(A,Hom_k(A,k^{(1)}))$$
and since for any algebra $B$ one has a fully faithful
embedding of $B-mod$ into $D^b(B)$ by considering a $B$-module as
a complex with differential $0$ and modules
concentrated in degree $0$ only,
this in turn gives a mapping
$$Hom_{D^b(A\otimes_kA^{op})}(A,Hom_k(A,k))\lra
Hom_{D^b(A\otimes_kA^{op})}(A,Hom_k(A,k^{(1)})).$$

Put $A^*:=Hom_k(A,k)$.
Recall that a $k$-algebra $A$ is symmetric if and only if there is
an isomorphism of $A\otimes_kA^{op}$-bimodules $A\simeq A^*$,
or equivalently there is a non degenerate symmetric bilinear form
$(\;,\;):A\times A\lra k$ satisfying $(a,cb)=(ac,b)$ for any
$a,b,c\in A$ (cf e.g. \cite[Chapter 9]{derbuch}).
Then the mapping $\zeta^A_n$
is defined by the equation
$(\zeta_n(z),x)^{p^n}=(z,x^{p^n})$ which can be written as
composition of the mappings in the following
diagram~$(\ddag)$:
$$
\begin{array}{cccccc}
&Hom_{D^b(A\otimes_kA^{op})}(A,A)&\lra&
Hom_{D^b(A\otimes_kA^{op})}(A,A^*) \\
&&& \downarrow
\mbox{\scriptsize $((Fr^k)_*)^n$}\\
(\ddag)\;\;\;&\uar\zeta_n^A&&
Hom_{D^b(A\otimes_kA^{op})}(A,Hom_{k}(A,k^{(n)})\\
&&&\uparrow
\mbox{\scriptsize $((Fr^A)^*)^n$}\\
&Hom_{D^b(A\otimes_kA^{op})}(A,A)&\lra&
Hom_{D^b(A\otimes_kA^{op})}(A,A^*)
\end{array}
$$
where the horizontal arrows are induced by the isomorphism
\begin{eqnarray*}
A&\lra&A^*\\
a&\mapsto&\left(b\mapsto (a,b)\right)
\end{eqnarray*}
which is coming from the symmetrising
bilinear form $(\phantom{x},\phantom{x}):A\otimes_kA\lra k$.
of $A$.

\section{Behaviour under derived equivalences}

\label{derivedinvariance}

In this section we prove our main result.

\begin{Theorem}\label{main}
Let $k$ be an algebraically closed field of characteristic $p>0$ and
let $A$ and $B$ be finite dimensional $k$-algebras. If
$D^b(A)\simeq D^b(B)$ as triangulated categories, then there is an
isomorphism $\varphi:ZA\lra ZB$ between the centres $ZA$ of $A$ and
$ZB$ of $B$ so that $\varphi(T_nA^\perp)=T_nB^\perp$ for all positive
integers $n\in\Z$.
\end{Theorem}

\begin{Rem}
This answers to the positive question 5.4 posed by
L\'aszl\'o H\'ethelyi, Ersz\'ebet Horv\'ath, Burkhard K\"ulshammer
and John Murray in \cite{Kuel}.
\end{Rem}

{\bf Proof:}
Let $F:D^b(A)\lra D^b(B)$ be a standard derived equivalence with
two-sided tilting complex $X$. Let $X'$ be the inverse tilting
complex. Then, in \cite{rogquest} it is shown that
$X\otimes_A-\otimes_AX'$ induces an equivalence
$G:D^b(A\otimes_kA^{op})\lra D^b(B\otimes_kB^{op})$ mapping the
$A$-$A$-bimodule $_AA_A$ to the $B$-$B$-bimodule $_BB_B$,
$$G(\ _AA_A)=\ _BB_B.$$
From \cite[Lemma 1]{rogquest} we know that
$$G\left(Hom_{k}(A,k)\right)=Hom_k(B,k).$$
We shall show that
$$G\left(Hom_{k}(A,k^{(n)})\right)=Hom_{k}(B,k^{(n)})$$
for all $n\in\Z$.
Indeed, $X\otimes_A-\simeq Hom_A(X',-)$ and
$-\otimes_AX'\simeq Hom_A(X,-)$ by the adjointness properties of
$Hom$ and $\otimes$-functors. Hence
(cf \cite[proof of Corollary 6.3.6]{derbuch}),
$$\begin{array}{lrcl}
&X\otimes_AHom_k(A,k^{(n)})\otimes_AX'&\simeq&
Hom_A(X',Hom_k(A,k^{(n)}))\otimes_AX'\\
&&\simeq&Hom_k(A\otimes_AX',k^{(n)})\otimes_AX'\\
&&\simeq&Hom_k(X',k^{(n)})\otimes_AX'\\
(\dagger)\;\;\;\;&&\simeq&Hom_A(X,Hom_k(X',k^{(n)}))\\
&&\simeq&Hom_k(X\otimes_AX',k^{(n)})\\
&&\simeq&Hom_k(B,k^{(n)})
\end{array}$$

We apply now $G$ to the diagram $(\ddag)$ of Section~\ref{sec1}
and get a commutative diagram
$$
\begin{array}{cccccc}
Hom_{D^b(B\otimes_kB^{op})}(B,B)&\lra&
Hom_{D^b(B\otimes_kB^{op})}(B,B^*)\\
{\mbox{\scriptsize $G$}}\uparrow \mbox{\scriptsize ${\simeq}$}&&
{\mbox{\scriptsize $G$}}\uparrow\mbox{\scriptsize ${\simeq}$}\\
Hom_{D^b(A\otimes_kA^{op})}(A,A)&\lra&
Hom_{D^b(A\otimes_kA^{op})}(A,A^*) \\
&& \downarrow
\mbox{\scriptsize $((Fr^k)_*)^n$}\\
\uar\zeta_n^A&&Hom_{D^b(A\otimes_kA^{op})}(A,Hom_{k}(A,k^{(n)}))\\
&&\uparrow
\mbox{\scriptsize $((Fr^A)^*)^n$}\\
Hom_{D^b(A\otimes_kA^{op})}(A,A)&\lra&
Hom_{D^b(A\otimes_kA^{op})}(A,A^*)\\
{\mbox{\scriptsize $G$}}\downarrow \mbox{\scriptsize ${\simeq}$}&&
{\mbox{\scriptsize $G$}}\downarrow\mbox{\scriptsize ${\simeq}$}\\
Hom_{D^b(B\otimes_kB^{op})}(B,B)&\lra&
Hom_{D^b(B\otimes_kB^{op})}(B,B^*)\;.
\end{array}
$$
It is clear that the upper and the lower square are commutative, since
they arise as squares induced from applying an equivalence of categories.

Recall the notation we use as explained at the
end of Section~\ref{somenotation}.

We obtain a commutative diagram
$$\begin{array}{ccc}
Hom_{D^b(A\otimes_kA^{op})}(A,Hom_k(A,k))&\stackrel{G}{\lra}&
Hom_{D^b(B\otimes_kB^{op})}(B,Hom_k(B,k))\\
\phantom{\mbox{\scriptsize $Hom_k(A,(Fr^k)^n)$}}\dar
\mbox{\scriptsize $Hom_k(A,(Fr^k)^n)$}&&\dar\varphi\\
Hom_{D^b(A\otimes_kA^{op})}(A,Hom_k(A,k^{(n)}))&\stackrel{G}{\lra}&
Hom_{D^b(B\otimes_kB^{op})}(B,Hom_k(B,k^{(n)}))\\
\end{array}$$
where $\varphi=G\circ Hom_k(A,(Fr^k)^n)\circ G^{-1}$.

We shall
need to see that $\varphi=Hom_k(B,(Fr^k)^n)$.

\begin{Claim}
$G\circ Hom(A,Fr^k) =Hom(B,Fr^k)\circ G$.
\end{Claim}

{\bf Proof:} Observe that $G=X\otimes_A-\otimes_AX'$ acts only on the
contravariant variables. Going through the isomorphisms
$(\dagger)$, since $Fr^k$ acts on the covariant variable only,
this proves the claim.
\hfill\dickebox

\medskip

Therefore, the diagram
$$
\begin{array}{cccc}
Hom_{D^b(B\otimes_kB^{op})}(B,B^*)&\stackrel{Hom(B,(Fr^k)^n)}{\lra}&
Hom_{D^b(B\otimes_kB^{op})}(B,Hom_{k}(B,k^{(n)} ))\\
\simeq\uar G&&\simeq\uar G\\
Hom_{D^b(A\otimes_kA^{op})}(A,A^*)&\stackrel{Hom(A,(Fr^k)^n)}{\lra}&
Hom_{D^b(A\otimes_kA^{op})}(A,Hom_{k}(A,k^{(n)} ))
\end{array}
$$
is commutative and the vertical morphisms are isomorphisms since $G$
is an equivalence, and since the images of the various objects under $G$
in their version $A$ and $B$ correspond to each other.

\begin{Claim}\label{postponed}
$G\circ Hom(Fr^A,k)=Hom(Fr^B,k)\circ G$.
\end{Claim}

Before starting with the proof observe the following consequences.
Once the claim is established the diagram
$$
\begin{array}{cccc}
Hom_{D^b(B\otimes_kB^{op})}(B,B^*)&\stackrel{Hom((Fr^B)^{n},k)}{\lra}&
Hom_{D^b(B\otimes_kB^{op})}(B,Hom_{k}(B,k^{(n)} ))\\
\simeq\uar G&&\simeq\uar G\\
Hom_{D^b(A\otimes_kA^{op})}(A,A^*)&\stackrel{Hom((Fr^A)^{n},k)}{\lra}&
Hom_{D^b(A\otimes_kA^{op})}(A,Hom_{k}(A,k^{(n)} ))
\end{array}
$$
is commutative.

Observe that since $Hom_k(Fr^A,k)$
is not $A\otimes_kA^{op}$-linear, the functor
$G$ is not defined on $Hom_k(Fr^A,k)$. Hence, the only way to prove the
commutativity of the above diagram is  by inspection
of the values.

\medskip

{\bf Proof of Claim~\ref{postponed}:}
We need to make explicit the mappings
$$G:Hom_{D^b(A\otimes_kA^{op})}(A,A^*)\lra
Hom_{D^b(B\otimes_kB^{op})}(B,B^*)$$
and
$$G:Hom_{D^b(A\otimes_kA^{op})}(A,Hom_{k}(A,k^{(1)} ))\lra
Hom_{D^b(B\otimes_kB^{op})}(B,Hom_{k}(B,k^{(1)})).$$
For this, it is useful, and possible, to replace $B$ by $X\otimes_AX'$ and
$A$ by $X'\otimes_BX$.

We first deal with the first identification.
Then, again by the usual adjointness formula between $Hom$ and $\otimes$,
one has to make explicit an isomorphism
$$G:Hom_{k}(A\otimes_{A\otimes_kA^{op}}A,k)\lra
Hom_{k}(B\otimes_{B\otimes_kB^{op}}B,k),$$
or, replacing $B$ by $X\otimes_AX'$ and $A$ by $X'\otimes_BX$,
$$G:Hom_{k}((X'\otimes_BX)\otimes_{A\otimes_kA^{op}}(X'\otimes_BX),k)\ra
Hom_{k}((X\otimes_AX')\otimes_{B\otimes_kB^{op}}(X\otimes_AX'),k).$$

The isomorphism
$$
A\otimes_{A\otimes_kA^{op}}A\simeq A/KA
\simeq B/KB\simeq B\otimes_{B\otimes_kB^{op}}B$$ comes from a
mapping $$(x\otimes y)\otimes 1_A\mapsto (y\otimes x)\otimes 1_B$$
where it is clear that this is well defined.
The diagonal mapping \begin{eqnarray*}
A&\lra&  \underbrace{
A\otimes_AA\otimes_A\dots\otimes_AA}_{p\mbox{ \scriptsize factors}}\\
a&\mapsto&a\otimes a\otimes\dots\otimes a
\end{eqnarray*}
is exactly the $p$-power map
$A\ni a\mapsto a^p\in A$. Composing with the natural projection $A\lra A/KA$
this defines the $p$-power mapping $A\ni a\mapsto a^p\in A/KA$.
If $k$ is of characteristic
$p$, then this last mapping is additive and factors through $A\lra A/KA$.

Now, we observe that
$A/KA$ is equally isomorphic to
$$\underbrace{A\otimes_AA\otimes_A\dots
\otimes_AA}_{p-1\mbox{ \scriptsize factors}}\otimes_{A\otimes_kA^{op}}A.$$
Moreover, we have seen $A\simeq X'\otimes_BX$ and
$B\simeq X\otimes_AX'$ and one recovers an isomorphism
\begin{eqnarray*}
(X'\otimes_BX)^{{\otimes_A}^{p-1}}\otimes_{A\otimes_kA^{op}}(X'\otimes_BX)&\lra&
(X\otimes_AX')^{{\otimes_B}^{p-1}}\otimes_{B\otimes_kB^{op}}(X\otimes_AX')\\
(x_i\otimes y_i)_{i=1}^{p-1}\otimes (x_p\otimes y_p)&\mapsto&
((y_p\otimes x_1)\otimes(y_{i}\otimes x_{i+1})_{i=1}^{p-2})
\otimes (y_{p-1}\otimes x_p)
\end{eqnarray*}
which is an incarnation of the isomorphism $A/KA\lra B/KB$.
We need to show that this is well-defined, but actually this is just a
straight forward and detailed examination which ring acts where in which way.

Therefore, the $p$-power map on $A/KA$ is mapped to the $p$-power map
on $B/KB$ by a standard derived equivalence.

We need to explain the second isomorphism
$$G:Hom_{D^b(A\otimes_kA^{op})}(A,Hom_{k}(A,k^{(1)} ))\lra
Hom_{D^b(B\otimes_kB^{op})}(B,Hom_{k}(B,k^{(1)} )).$$
Here, we observe that
$$Hom_{D^b(A\otimes_kA^{op})}(A,Hom_{k}(A,k^{(1)} ))\simeq
Hom_{k}(A\otimes_{A\otimes_kA^{op}}A,k^{(1)})$$ and the very same arguments
and constructions as above hold. The only difference is that one needs to
consider semilinear mappings only at the end. The reorganization
procedure is just the same.
In particular, the action of $Fr^B$ consists
in tensoring the whole term on the right $p$ times over $B\otimes B^{op}$.
It is now immediate to see that
this operation commutes with this reorganization of factors
as described by explaining $\nu$.
So, $$G\circ Hom(Fr^A,k)=Hom(Fr^B,k)\circ G.$$
\hfill\dickebox

\begin{Claim}
The images of $G\circ \zeta_n^A\circ G^{-1}$ and of
$\zeta_n^B$ coincide.
\end{Claim}

{\bf Proof:}
Since $\varphi:A\lra Hom_k(A,k)$ is an isomorphism
of $A\otimes_kA^{op}$-modules, and since $G$ is a functor, $G\varphi$
is an isomorphism as well.
As we know that choosing an isomorphism $B\lra Hom_k(B,k)$ is equivalent
to choosing a symmetrising form making $B$ into a symmetric algebra,
we may well work with this form instead of the original one. Actually,
given two different isomorphisms $\phi:B\lra Hom_k(B,k)$ and $\psi:B\lra
Hom_k(B,k)$, then for all $x\in B$ one has
$\phi^{-1}\psi(x)=\lambda x$ for an invertible central
$\lambda\in Z(B)^{*}$.  So, the resulting $\zeta^B_n$ differ by
invertible central elements. As a consequence, the images are
identical.
\hfill\dickebox

\medskip

We shall finish the proof of the theorem.
By the previous claims the composition
\begin{eqnarray*}
ZB\stackrel{G\phi}{\lra}Hom_k(B/KB,k)&
\stackrel{(G(Fr^A)^*G^{-1})^n}{\lra}&
Hom_k(B/KB,k^{(n)})\\
&\stackrel{((G(Fr^k)^*)^n)^{-1}}{\lra}&Hom_k(B/KB,k)
\stackrel{G\phi^{-1}}{\lra}ZB
\end{eqnarray*}
is a mapping which differs from $\zeta_n^B$ by some central unit of
$B$ and therefore the isomorphism induced by $G$ between the
centres of $A$ and $B$
maps $T_nA^{\perp}$ to $T_nB^{\perp}$.

This finishes the proof of the theorem.\hfill\dickebox

\end{document}